\documentclass[9pt,leqno]{amsart}

 \usepackage{mathrsfs}

\usepackage{color}
\usepackage{amsmath, amsthm, amssymb}
\usepackage{amsfonts}
\usepackage[dvips]{epsfig}
\usepackage{graphicx}
\usepackage{caption}
\usepackage{subcaption}
\usepackage[english]{babel}
\usepackage[left=4.5cm,right=4.5cm,top=2cm,bottom=2cm]{geometry}
\usepackage{hyperref}

\usepackage{tikz}
\usepackage{rotating}
\usepackage[utf8]{inputenc}
\usepackage{cite}
\usepackage{amscd}
\usepackage{color}
\usepackage{bm}
\usepackage{enumerate}

\usepackage{verbatim}
\usepackage{hyperref}
\usepackage{amstext}
\usepackage{latexsym}
\theoremstyle{plain}
\newtheorem{theorem}{Theorem}[section]

\newtheorem{corollary}[theorem]{Corollary}
\newtheorem{proposition}[theorem]{Proposition}
\newtheorem{lemma}[theorem]{Lemma}
\newtheorem{remark}[theorem]{Remark}

\numberwithin{theorem}{section}
\numberwithin{equation}{section}

\def\square{\vbox{
    \hrule height .4pt
    \hbox{\vrule width .4pt height 7pt \kern 7pt
       \vrule width .4pt}
    \hrule height .4pt }}

\newcommand{\average}{{\mathchoice {\kern1ex\vcenter{\hrule height.4pt
width 6pt depth0pt} \kern-9.7pt} {\kern1ex\vcenter{\hrule
height.4pt width 4.3pt depth0pt} \kern-7pt} {} {} }}

\def\R{\mathbb{R}}

\renewcommand{\d}{\delta }

\newcommand{\D }{\Delta }

\newcommand{\e }{\varepsilon }

\renewcommand{\l }{\lambda }

\newcommand{\n }{\nabla }
\newcommand{\vp }{\varphi }
\renewcommand{\phi}{\varphi}

\newcommand{\s }{\sigma }

\renewcommand{\O }{\Omega }

\newcommand{\be}{\begin{equation}}
\newcommand{\ee}{\end{equation}}

\newcommand{\de}{\partial}

\newcommand{\calC }{\mathcal{C}}

\newcommand{\calD }{\mathcal{D}}

\newcommand{\B}{{Q}}

\renewcommand{\epsilon}{\varepsilon}

\begin{document}
\author{Mamadou Ciss}
\address{M. C. : Université Alioune Diop, Bambey.}
\email{mamadou.ciss@uadb.edu.sn}
\author{Abdourahmane Diatta}
\address{A.D. : Université Assane Seck de Ziguinchor, UFR des Sciences et Technologies, département de mathématiques, Ziguinchor.}
\email{a.diatta20160578@zig.univ.sn}
\author{El Hadji Abdoulaye THIAM}
\address{H. E. A. T. : Université Iba Der Thiam de Thies, UFR des Sciences et Techniques, département de mathématiques, Thies.}
\email{elhadjiabdoulaye.thiam@univ-thies.sn}
\title[A Nonlinear elliptic PDE with curve singularity on the boundary]
{A Nonlinear elliptic PDE with curve singularity on the boundary}

\begin{abstract}
Let $\Omega$ be a bounded domain of $\mathbb{R}^{N+1}$ ($N \geq 3$) with smooth boundary $\partial \Omega$ and $\Sigma$ be a closed submanifold contained on $\partial \Omega$ and containing $0$. We are interesting in the existence of positive $H^1(\Omega)$-solution of the following Hardy-Sobolev trace type equation
\begin{equation*}
\begin{cases}
-\Delta u+u=0 \qquad & \textrm{ in $\Omega$}\\\\
\displaystyle\frac{\partial u}{\partial \nu}= \rho_{\Sigma}^{-s} u^{q_s-1} \qquad & \textrm{ on $\partial \Omega$},
\end{cases}
\end{equation*}
where $\nu$ is the unit outer normal of $\partial \Omega$, $\rho_\Sigma: \partial \Omega \to \mathbb{R}$ is the distance function in $\partial \Omega$ to the curve $\Sigma$:
$$
\rho_\Sigma(x):= \inf_{y \in \Sigma} d_{\tilde{g}}(x, y)
$$ and for $0\leq s <1$, $q_s:=\frac{2(N-s)}{N-1}$ is the critical Hardy-Sobolev exponent. The existence of solution may depend on the local geometry of the boundary $\partial \Omega$ and $\Sigma$ at $0$ or in the shapes of the domain $\Omega$ and its boundary $\partial \Omega$.
\end{abstract}
\maketitle


%
\section{Introduction and Main Results}
Let $\O$ be a bounded domain of $\mathbb{R}^{N+1}$ with $N \geq 2$ with boundary $\partial \Omega$ and $\Sigma$ be a closed regular curve contained on $\partial \Omega$. We assume that $0\in \Sigma$. We consider $(\partial \Omega, \tilde{g})$ as a
Riemannian manifold, with metric $\tilde{g}$ induced by $\mathbb{R}^{N+1}$ on the boundary $\partial \Omega$. Given $s\in [0, 1)$, we study existence of positive $H^1(\Omega)$-solution for the following Hardy-Sobolev trace equation with singularity a curve
\begin{equation}\label{Hardy-Sobolev-Trace-Inequality}
\begin{cases}
-\Delta u+u=0 \qquad & \textrm{ in $\Omega$}\\\\
\displaystyle\frac{\partial u}{\partial \nu}= \rho_{\Sigma}^{-s} u^{q_s-1} \qquad & \textrm{ on $\partial \Omega$}
\end{cases}
\end{equation}
where $\nu$ is the unit outer normal of $\partial \Omega$, $q_s=\frac{2(N-s)}{N-1}$ is the critical Hardy-Sobolev exponent and  $\rho_\Sigma: \partial \Omega \to \mathbb{R}$ is the distance function in $\partial \Omega$ to the curve $\Sigma$. That is
\begin{equation}\label{Def-Dist-On-Boundary}
\rho_\Sigma(x):=d_{\tilde{g}}(x, \Sigma):= \inf_{y \in \Sigma} d_{\tilde{g}}(x, y).
\end{equation}
Solutions of \eqref{Hardy-Sobolev-Trace-Inequality} are minimizers for the functional $J: H^1(\Omega) \setminus \lbrace 0\rbrace \to \R$ defined by
$$
J(u):= 
\frac{\displaystyle\int_{\Omega} \left(|\nabla u|^2+ u^2\right) dx}{\left(\displaystyle\int_{\partial \Omega} \rho^{-s}_{\Sigma} |u|^{q_s} d\sigma(x)\right)},
$$
where $H^1(\Omega)$ is the completion of $\mathcal{C}^\infty_c(\Omega)$ with respect to the norm
$$
u \longmapsto \sqrt{\displaystyle\int_{\Omega} \left(|\nabla u|^2+ u^2\right) dx}.
$$
Next, we set
$$
S_{\Omega, \Sigma, s}:=\inf_{u \in H^1(\Omega) \setminus \lbrace 0 \rbrace} J(u).
$$
The exponent $q_s$ defined above is critical in the sens that $H^1(\Omega)$ is continuously embedded in
$$
L^p(\Omega, \rho^{-s}_{\Sigma})
:=
\lbrace
u: \Omega \to \mathbb{R}: \quad \textrm{ such that} \int_{\partial \Omega} \rho^{-s}_{\Sigma} |u|^{p} d\sigma(x) <\infty
\rbrace
$$
if $1 \leq p \leq q_s$ and the embedding is compact if $1 \leq p <q_s$. We refer  to Ghoussoub-Kang \cite{Ghoussoub-Kang} for more details. Then thanks to the above embedding, the functional $J$ is well defined. 
In this kind of problem, the difficulty is the lack of compactness, for the critical embedding. To recover compactness, we will prove that there exists $S_{N, s}$ positive depending on $N$ and $s$ such that if
\begin{equation}\label{Compactness-Ineq}
S_{\Omega, \Sigma, s}<S_{N, s}
\end{equation}
then any minimizing sequence $(u_n) \subset H^1(\Omega)$ admits a subsequence  that converges strongly in $H^1(\Omega)$ to some function $u \neq 0$. The constant $S_{N, s}$ is related to the limiting Hardy-Sobolev trace equation $\mathbb{R}^{N+1}_+$. Indeed,  we let $x=(t, y, z) \in \mathbb{R} \times \mathbb{R}^{N-1} \times \mathbb{R}_+$  and  consider the Euler-Lagrange equation
$$
\begin{cases}
\Delta w=0 \qquad &\textrm{ in $\mathbb{R}\times \mathbb{R}^{N-1} \times \mathbb{R}_+$}\\\\
\displaystyle-\frac{\partial w}{\partial z}=S_{N, s} w^{q_s-1} |y|^{-s}  &\textrm{ on  $\mathbb{R}\times \mathbb{R}^{N-1}$}.
\end{cases}
$$
Defining $\Phi: \mathcal{D}^{1,2}(\mathbb{R}^{N+1}_+) \to \mathbb{R}$ by:
$$
\Phi(u)=\frac{\displaystyle \int_{\mathbb{R}^{N+1}_+} |\nabla u(t, y, z)|^2 dt dy dz}{\displaystyle \left(\int_{\mathbb{R}^{N}}|y|^{-s} |u|^{q_s} dt dy \right)^{2/q_s}}.
$$
Then the constant $S_{N, s}$ defined in \eqref{Compactness-Ineq} is given by
\begin{equation}\label{HS-Trace}
 S_{N, s}:= \inf_{u \in \mathcal{D}^{1, 2}(\mathbb{R}^{N+1}_+)} \Phi(u)
\end{equation}
where we recall that $\mathbb{R}^{N+1}_+=\lbrace x=(t, y, z)\in \mathbb{R} \times \mathbb{R}^{N-1} \times \mathbb{R} \rbrace$ and $\mathcal{D}^{1, 2}(\mathbb{R}^{N+1}_+$ is the completion of $\mathcal{C}^\infty_c(\mathbb{R}^{N+1}_+)$ with respect to the norm
$$
\displaystyle u \longmapsto \sqrt{\int_{\mathbb{R}^{N+1}_+} |\nabla u|^2 d x}.
$$
In the sequel, we define by $H_{\partial \Omega}$ the mean curvature of the boundary $\partial \Omega$ and by $\mathcal{H}_{\partial \Omega}$ and $\mathcal{H}_1$ two others geometric quantities depending on the principal curvatures, see for instance Section \ref{Section2} below. Then as a consequence of inequality  \eqref{Compactness-Ineq}, we have our first main result where the existence of solution depends on the local geometries of the boundary $\partial \Omega$ and of the curve $\Sigma$.
\begin{theorem}\label{Theorem1}
We consider a smooth bounded domain
 $\O$ of $\mathbb{R}^{N+1}$ with $N\geq 3$, $\Sigma$ a smooth closed curve contained in $\partial \Omega$  with $0\in\Sigma$ and $s\in[0,1)$. Assume that
	\begin{equation}\label{Geometrique}
	H_{\partial \Omega}(0)-A_{N, s}\mathcal{H}_{\partial \Omega}(0)-B_{N, s}\mathcal{H}_1  <0.
	\end{equation} 
	with
	$$
	A_{N, s}:=2\frac{\displaystyle\int_{\mathbb{R}^{N+1}_+} z |\nabla_y w|^2 dx}{\displaystyle\int_{\mathbb{R}^{N+1}_+} z |\nabla w|^2 dx} \qquad \textrm{ and } \qquad B_{N, s}:= \frac{\displaystyle\int_{\mathbb{R}^{N+1}_+} z |\frac{\partial w}{\partial t}|^2 dx}{\displaystyle\int_{\mathbb{R}^{N+1}_+} z |\nabla w|^2 dx}.
	$$
	Then $S_{\Omega, \Sigma,s}$ is achieved by a positive function $u\in H^1(\Omega)$ satisfying
	\begin{align*}
	\begin{cases}
	\displaystyle-\Delta u+u=0&\qquad \textrm{ in } \Omega \vspace{3mm}\\
	\displaystyle \frac{\partial u}{\partial \nu}=S_{\Omega, \Sigma, s}\, \rho_\Sigma^{-s}\,  u^{q_s-1}& \qquad \textrm{ on } \partial \Omega.
	\end{cases}
	\end{align*}
\end{theorem}
Next, we turn out to our second main result where the existence of solution does not depend on the local geometry of the domain but on the shape of the domain $\Omega$  and its boundary $|\partial \Omega|$. More precisely, we have
\begin{theorem}\label{Theorem2}
	Let $\Omega$ be a smooth bounded domain of $\mathbb{R}^{N+1}$, $N\geq 2$, $\Sigma$ a smooth closed curve contained in $\partial \Omega$  with $0\in\Sigma$ and $s\in[0,1)$. Assume that
	\begin{equation}\label{Algebrique}
	\frac{\displaystyle |\Omega|}{\left( \displaystyle \int_{\partial \Omega} \rho_\Sigma^{-s} d\sigma(x) \right)^{2/q_s}} <S_{N, s},
	\end{equation}
	where $|\Omega|$ is the volume of $\Omega$.
	Then $\mathcal{S}_{\Omega, \Sigma, s}$ is achieved by a positive function $u\in H^1(\Omega)$ satisfying
	\begin{align*}
	\begin{cases}
	\displaystyle-\Delta u+u=0&\qquad \textrm{ in } \Omega \vspace{3mm}\\
	\displaystyle \frac{\partial u}{\partial \nu}=S_{\Omega, \Sigma, s}\, \rho_\Sigma^{-s}\,  u^{q_s-1}& \qquad \textrm{ on } \partial \Omega.
	\end{cases}
	\end{align*}
\end{theorem}

The proofs of Theorems \ref{Theorem1} and \ref{Theorem2} are based  on test function method. We construct suitable functions that allow us to compare $S_{\O, \Sigma, s}$ with $S_{N, s}$. Since the inequality $S_{\Omega, \Sigma, s} \leq S_{N, s}$ always holds, the main goal is to identify test functions such that, under the assumptions of Theorems\ref{Theorem1} and \ref{Theorem2}, relation \eqref{Compactness-Ineq} is satisfied.  This enables us to recover compactness, ensuring that every minimizing sequence for $S_{\O, \Sigma, s}$ possesses a subsequence converging to a minimizer.\\

This paper proceeds as follows. In Section \ref{Section2}, we introduce Fermi coordinates, derive the expansion of the metric, and compute its inverse and determinant. Section \ref{Section3} deals with the existence, symmetry, and other properties of the ground-state solution of the Hardy-Sobolev equation in the Euclidean setting. In Section \ref{Section4}, we prove existence of minimizer for the Hardy-Sobolev trace best constant in domains and in Section \ref{Section5} we build the appropriate test functions  and conclude with the proofs of the main results.

\section{Geometric Preliminaries}\label{s:Geometric-prem}\label{Section2}
\subsection{Fermi-coordinates}
Let $\Omega$ be a bounded domain of $\R^{N+1}$, with $N \geq 2$ and let $\Sigma$ be a regular closed curve contained in $\de \Omega$, the boundary of $\Omega$. We consider $(\de \O, \tilde{g}))$ as a
Riemannian manifold, with Riemannian metric $\tilde{g}$ induced by $\R^{N+1}$ on $\de \O$ and denote by $\rho_\Sigma: \de \O \to \R_+$  the Riemannian distance in $(\de \O, \tilde{g}))$ to the curve $\Sigma$:
$$
\rho_\Sigma(x):=\inf_{y \in \Sigma} d_{\tilde{g}}(x, y).
$$
We assume that $0\in \Sigma$ and we have the natural splliting
$$
T_0\de \Omega=T_0 \Sigma\oplus N_0 \Sigma,
$$
where $T_0\de \Omega$ is the tangent space of $\de \O$ at $0$, $T_0\Sigma$ is the tangent space of $\Sigma$ and $N_0\Sigma$ is the normal bundle of $\Sigma$. Next, we choose orthonormal bases
$(E_0)$ of $T_0\Sigma$ and $(E_2, \cdots, E_N)$ of $N_0\Sigma$.

For $r>0$ small enough, we consider the regular curve $\gamma:(-r, r)\to \Sigma$ such that $\gamma(0)=0_{\R^{N+1}}$. This yields the coordinate vector field
$$
X_1(t):=\gamma_*\left(\frac{\de }{\de t}\right),
$$
so that
\begin{equation}\label{Ami0}
\nabla_{X_1} X_1|_{0} \in N_0 \Sigma.
\end{equation}
 Therefore, there exists real numbers $(\beta_i)_{1\leq i \leq N}$ such that
\begin{equation}\label{Ami1}
\nabla_{X_1} X_1|_{0}= \sum_{i=2}^N \beta_i E_i.
\end{equation}
The vectors $(E_i)_{2 \leq i \leq N}$ are extended along the curve $\gamma(t)$ by parallel transport with respect to the induced connection on the normal bundle $N\Sigma$, thereby defining an orthonormal frame field $(X_i)_{2 \leq i \leq N}$ of $N\Sigma$ in a neighborhood of $0$ in $\Sigma$ which satisfies 
\begin{equation}\label{Ami2}
\n_{X_1} X_i\mid_{0} \in T_0 \Sigma
\end{equation}
and hence for $i=2, \cdots, N$, there exists $\kappa_i$ real number satisfying
\begin{equation}\label{Ami3}
\n_{X_1} X_i\mid_{0}=\kappa_i E_1.
\end{equation}
 We introduce geodesic normal coordinates in a neighborhood of $0 \in \Sigma \subset \de \O$ with coordinates $(t, y):=(t, y_2, \cdots, y_N) \in \R^{N}$. We set
$$
f(t, y)=\textrm{Exp}_{\gamma(t)}^{\de \O} \left(\sum_{i=2}^N y_i X_i(t)  \right), \qquad (t, y)\in (-r, r)\times B_{\R^{N-1}}(0, r),
$$
where $B_{\R^{N-1}}(0, r)$ is the ball of $\R^N$ centered at the origin and of radius $r>0$. It is clear that in this geodesic local coordinates, the distance function $\rho_\Sigma$ satisfies
\begin{equation}\label{Heat1}
\rho_\Sigma(f(t, y))=|y|.
\end{equation}
\subsection{The Induced Metric $\tilde{g}$ on the boundary $\de \O$}
From this choice of coordinates, we obtain the corresponding coordinate vector fields on $\partial \Omega$:
$$
X_{1}=f_*\left(\frac{\de }{\de t}\right) \quad\textrm{and} \quad    X_i(t, y)=f_*\left(\frac{\de }{\de y_i}\right) \quad \textrm{ for $i=2, \cdots, N$}.
$$ 
For $i, j=1, \cdots, N$, we let $\tilde{g}_{ij}(x)=\langle X_i , X_j \rangle$  to be the components of the metric $\tilde{g}$ on $\de \O$. By construction, we have
$$
X_i|_0=E_i \qquad \textrm{ for all $i=1, \cdots, N$}.
$$ Then near the origin we have the
\begin{proposition}\label{Metric-Induced}
Let $\tilde{x}:=(t, y) \in Q_r$ and $i, j=2, \cdots, N$. Then we have
$$
\begin{cases}
\tilde{g}_{11}(\tilde{x})&=\displaystyle 1+ 2\sum_{i=2}^{N} \kappa_i y_i + O(|(\tilde{x})|^2)\\\\
\tilde{g}_{i1}(\tilde{x})&= \beta_i t+O(|\tilde{x}|^2)\\\\
\tilde{g}_{ij}(\tilde{x})&= \displaystyle\delta_{ij}   + O(|\tilde{x}|^2),
\end{cases}
$$
where we have denoted 
\begin{equation}\label{Cube-Cube}
Q_r:=(-r, r) \times B_{\R^{N-1}}(0, r).
\end{equation}
\end{proposition}

\begin{proof}
The proof of the result is divided into several steps. Each step corresponds to a class of components of the metric.

\textbf{Step 1: Computation of component $\tilde{g}_{11}$.} By Taylor expansion (around the origin), we have
\begin{align}\label{Step1-0}
&\displaystyle g_{11}(t, y)= \left\langle X_1, X_1 \right\rangle|_0 + X_1\left\langle X_1, X_1\right\rangle |_0 t+ \sum_{i=2}^{N}X_i\left\langle X_0, X_0\right\rangle|_0 y_i+O(|\tilde{x}|^2).
\end{align}
Note that $\left\langle X_1, X_1 \right\rangle|_0=\langle E_1, E_1 \rangle=1$. Moreover by \eqref{Ami0}, we have
\begin{align}\label{Step1-1}
X_1\left\langle X_1, X_1 \right\rangle|_0= 2\left\langle \n_{X_1} X_1, X_1  \right\rangle |_0=0.
\end{align}
By \eqref{Ami3}, we have, for $i=2,\cdots,N$, that 
\begin{align}\label{Step1-20}
X_i\left\langle X_1, X_1\right\rangle &= 2\left\langle \n_{X_i} X_1, X_1\right\rangle|_0=2 \kappa_i.
\end{align} 
Therefore by \eqref{Step1-0}, \eqref{Step1-1} and \eqref{Step1-20}, we get
\begin{align}\label{Step1-2}
\displaystyle g_{11}(\tilde{x})&= 1+ 2 \sum_{i=2}^N \kappa_i y_i+O(|\tilde{x}|^2).
\end{align}

\textbf{Step 2: Computation of component $\tilde{g}_{1i}$}. In this part, we assume that $i=2,\cdots,N$. Then we have 
	\begin{align*}
	\displaystyle
	\tilde{g}_{i1}(x)&= \langle X_1, X_i \rangle|_0 +X_1\left\langle X_i, X_1\right\rangle t+\sum_{j=2}^N X_j\left\langle X_i, X_1\right\rangle y_j+O(|\tilde{x}|^2).
	\end{align*}
	Note that $\langle X_1, X_i \rangle|_0=0$ and by  \eqref{Ami1} and \eqref{Ami3}, we have 
	\begin{align*}
	X_1\left\langle X_i,X_1\right\rangle&= \displaystyle \left\langle \n_{X_1}X_i, X_1\right\rangle +\left\langle X_i,\n_{X_1} X_1\right\rangle= \beta_i.
	\end{align*}
Next, for $i, j=2, \cdots, N$, we get
	\begin{align*}
	X_j\left\langle X_i, X_1\right\rangle& =  \left\langle \n_{X_j} X_i,X_1\right\rangle + \left\langle X_i,\n_{X_j} X_1 \right\rangle=0
	\end{align*}
	thanks to \eqref{Ami3}.
Therefore 
$$
\tilde{g}_{i1}= \beta_i t+o(|\tilde{x}|^2).
$$

\textbf{Step 3: Computations of the components $\tilde{g}_{ij}$}. For $i, j=2, \cdots, N$, we have	
	 \begin{align*}
	  	\displaystyle \tilde{g}_{ij}(x)= \delta_{ij} +X_1\left\langle X_i,X_j\right\rangle t+ \sum_{k=2}^NX_k\left\langle X_i,X_j\right\rangle y_k  +O(|\tilde{x}|^2). 
	 \end{align*}
By \eqref{Ami3} , we have 
	 \begin{align*}
	 X_1 \left\langle X_i,X_j\right\rangle= \left\langle\n_{X_1} X_i,X_j\right\rangle + \left\langle X_i, \n_{X_1}X_j\right\rangle= 0
	 \end{align*}
and for $k=2, \cdots,\cdots, N$, we have 
	 \begin{align*}
	 X_k\left\langle X_i,X_j\right\rangle =
	 \left\langle \n_{X_k}X_i,X_j\right\rangle + \left\langle X_i,\n_{X_k} X_j\right\rangle=0.
	 \end{align*}
	Therefore
	$$
	\tilde{g}_{ij}= \delta_{ij} +O(\tilde{x}|^2).
	$$
	This then ends the proof.	
\end{proof}
\subsection{The Full Metric}
Hereafter, $N_{\Sigma}$ denotes the unit normal vector field along the curve, directed toward the interior of $\Omega$. Up to rotation, we may assume that $N_{\Sigma}(0)=e_{N+1}$. For any vector field $X$ on $T \de \O$, we define
$$
H(X)=dN_{\de \O}[X].
$$
Then the "normalized" \textbf{mean curvature} of $\Sigma$ at $0$ is given by
\begin{equation}\label{Mean-Curvature}
H_{\de \Omega}(0)=\frac{1}{N}\sum_{i=1}^{N}  \langle H(E_i), E_i \rangle. 
\end{equation}
In the sequel, we denote by $\mathcal{H}_{\de \O}(0)$ and $\mathcal{H}_1$ the following geometric quantities
\begin{equation}\label{Rel-Mean-Curvature}
\mathcal{H}_{\de \Omega}(0)=\frac{1}{N-1}\sum_{i=2}^{N}  \langle H(E_i), E_i \rangle \quad \textrm{ and } \quad \mathcal{H}_1(0)= \langle H(E_1), E_1 \rangle.
\end{equation}
Now we consider a local parametrization of a neighborhood of $0$ in $\R^{N+1}$ defined as
\begin{equation}\label{Full-Parametrization}
F(t, y, z)=f(t, y)+ z N_{\Sigma} \left(f(t, y)\right), \qquad \forall (t, y, z) \in \mathcal{Q}_r:=Q_r\times (0, r),
\end{equation}
where $Q_r$ is defined in \eqref{Cube-Cube}.
This yields the coordinates vector-fields in $\R^{N+1}$ as:
$$ 
Y_1(t, y, z)=F_*\left(\frac{\de }{\de t}\right);\quad Y_{N+1}(t, y, z)=F_*\left(\frac{\de }{\de z}\right)    \quad \textrm{and} \quad Y_i(t, y, z)=F_*\left(\frac{\de }{\de y_i}\right) \textrm{ for $i=2, \cdots, N$}.
$$
Then near the point $F(t, y, 0)$, we have
\begin{equation}\label{Last-Component}
Y_i=X_i+z H(X_i)+O(|z|^2),
\end{equation}
see for instance [\cite{Fall}, Lemma 2.1]. 
Let 
$$
g_{ij}(x):= \langle Y_i, Y_j\rangle \qquad \textrm{ for $i, j=1, \cdots, N+1$}
$$
denote the components of the flat Riemannian metric $g$. Then we have the following expansion.

\begin{proposition}\label{Full-Metriquea}
Let $x:=(t, y, z) \in \mathcal{Q}_r$ and $i, j=2, \cdots, N$. Then we have
$$
\begin{cases}
\displaystyle g_{11}(x)=1+ 2\sum_{i=2}^{N} \kappa_i y_i +2z \langle H(X_1), X_1\rangle+O(|x|^2)\\\\
g_{i1}(x)=\beta_i t+2z\langle H(X_1, X_i\rangle+O(|x|^2)\\\\
g_{i N+1}(x)= O(|x|^2)\\\\
g_{ij}(x)=\delta_{ij}+2 z\langle H(X_i, X_j\rangle+O(|x|^2)\\\\
g_{1 N+1}(x)=g_{iN+1}(x)=O(|x|^2)\\\\
\displaystyle g_{N+1 N+1}(x)=1.
\end{cases}
$$
\end{proposition}
\begin{proof}
The proof is an immediate consequence of Proposition \ref{Metric-Induced} and identity \eqref{Last-Component}.
\end{proof}
\subsection{The inverse and the determinant of the Riemannian metric $g$}
As a first consequence of Proposition \ref{Full-Metriquea}, rewritten the matrix $g$ in the form
\begin{equation}\label{B1}
g=I_{N+1}-A,
\end{equation}
where $I_{N+1}$ is the identity matrix  of $\mathcal{M}_{N+1}(\R)$. Since $|A|=O(r)$, then for $r$ small enough $g$ is invertible. Moreover we have
\begin{corollary}\label{Inverse-Metric-Cor1}
Let $x:=(t, y, z) \in \mathcal{Q}_r$ and $i, j=2, \cdots, N$. Then the components $\left(g^{\alpha \beta} \right)_{1\leq \alpha, \beta \leq N+1}$ of the inverse metric $g$ are given by
$$
\begin{cases}
\displaystyle g^{11}(x)=1- 2\sum_{i=2}^{N} \kappa_i y_i -2z \langle H(X_1), X_1\rangle+O(|x|^2)\\\\
g^{i1}=-\beta_i t-2z\langle H(X_1, X_i\rangle+O(|x|^2)\\\\
g^{i N+1}(x)= O(|x|^2)\\\\
g^{ij}=\delta_{ij}-2 z\langle H(X_i, X_j\rangle+O(|x|^2)\\\\
g^{1 N+1}(x)=g_{i N+1}(x)=O(|x|^2)\\\\
\displaystyle g_{N+1 N+1}(x)=1.
\end{cases}
$$
\end{corollary}
\begin{proof}
By \eqref{B1}, we can deduce that, for $r$ small enough, the matrix $g$  is invertible. Moreover its inverse is given by the power series
$$
g^{-1}=\sum_{k=0}^{\infty} A^k= I_{N+1}+A+O(|A|^2).
$$
Consequently, the result is an immediate consequence of this identity together with Proposition~\ref{Full-Metriquea}. This completes the proof.
\end{proof}
We finish this section by given the Taylor expansion of the root of the determiant of $g$ for $r$ small enough. We have
\begin{corollary}\label{Det-Metrci-Cor2}
Let $x:=(t, y, z) \in \mathcal{Q}_r$. Then for $r$ small enough, the determiant $|g|$ of the metric $g$ satisfies:
$$
\sqrt{|g|}(x)=1+\sum_{i=2}^{N} \kappa_i y_i+z H_{\de \O}(0)+O(|x|^2).
$$
\end{corollary}
\begin{proof}
For $r$ small enough, we have the following expansion
\begin{equation}\label{Yaar}
\sqrt{|g|}(x)=1+\frac{tr(A)}{2}+O(|A|^2),
\end{equation}
as $|A| \to 0$, where 
$$
tr(A)=\sum_{i=1}^{N+1} A_{i i }
$$ is the trace of the matrix $A$. Then we get the result immediately from \eqref{Yaar} and Proposition \ref{Full-Metriquea}.
\end{proof}
\section{The limiting problem}\label{Section3}
\subsection{Existence of ground states in $\R\times \R^{N-1}\times \R_+$}\label{ss:ex-gs}

\begin{theorem}\label{th:exist-s01}
Let $N \geq 2$, $x=(t, y, z) \in \mathbb{R}^{N+1}_+$,  $s\in(0,1)$.  Consequently, there exists $w\in \calD:=\calD^{1,2}(\R^{N+1}_+)$ solution of
\begin{align*}
\begin{cases}
\displaystyle \D w =0&\qquad \textrm{ in } \mathbb{R}^{N+1}_+, \vspace{3mm}\\
\displaystyle -\frac{\de w}{\de z}=S_{N, s} w^{q_s-1}|y|^{-s} &\qquad \textrm{ on } \R^N, \vspace{3mm}\\
\displaystyle\int_{\R^{N}} |y|^{-s} w^{q_s} dy dt=1.
\end{cases}
\end{align*}

 \end{theorem}
\begin{proof}
We start by defining the functionals	 $\Phi,\Psi: \calD\to \R$ by
$$
\Phi(w):=\frac{1}{2} \int_{\R^{N+1}_+}|\n w|^2dx  \qquad \textrm{and} \qquad
\Psi(w)=\frac{1}{q_s}\int_{\R^{N}}|y|^{-s} |w|^{q_s}dy dt.
$$
Using Ekland variational principle, we can easily prove the existence of a minimising sequence  $(w_n)_n$  for $ S_{N,s}$. In others words, the sequence $(w_n)_n$ satisfies
\begin{equation}\label{eq:wenorm}
\int_{\R^{N}}|y|^{-s} |w_n|^{q_s}d y dt=1 \qquad \textrm{and} \qquad
\Phi(w_n)=\frac{1}{2}S_{N, s}+o_n(1).
\end{equation}
Moreover, we have
\begin{equation}\label{eq:weps-stf}
\Phi'(w_n)=S_{N,s}\Psi'(w_n)++o_n(1)\qquad\textrm{ in } \calD',
\end{equation}
where $o_n(1) \to 0$ as $n \to \infty$ and $\calD' $ denotes the dual space of  $\calD$.
By the second identity in  \eqref{eq:wenorm}, there exists $C>0$ such that
\begin{equation}\label{eq:webndH1}
  \int_{ \R^{N+1}_+}|\n w_n|^2dz  \leq C.
\end{equation}
Let $Q:(0,\infty)\to\mathbb{R}$ denote the Levi-type concentration function, given by
\[
Q(r) := \int_{B^N_r} |y|^{-s} |w_n|^{q_s} \, dy \, dt.
\]
Using the continuity of $Q$ and \eqref{eq:wenorm}, there exists $r_n>0$ such that
$$
Q(r_n):=\int_{ B^N_{r_n}}|y|^{-s} |w_n|^{q_s}d y dt=\frac{1}{2}.
$$
Next, we set $v_n(x):=r_n^{\frac{N-1}{2}}w_n(r_n x)$. Then, we have
$$
\int_{ \R^{N+1}_+}|\n w_n|^2dx=\int_{ \R^{N+1}_+}|\n v_n|^2dx,\quad \int_{ \R^N}|y|^{-s} |w_n|^{q_s}d ydt=\int_{ \R^N}|y|^{-s} |v_n|^{q_s}dy dt 
$$
and
\begin{equation}\label{eq:Levi}
\int_{ B^N_1}|y|^{-s} |v_n|^{q_s}d y dt=\frac{1}{2}.
\end{equation}
Consequently $(v_n)_n$ is a minimizing sequence. In particular, there exists $v \in \calD$ such that   ${v_n} \rightharpoonup v$. 
In the sequel, we will show that $v\neq0$. Otherwise, we have $v_n\to0$ in $L^2_{loc}(\R^{N+1}_+)$ and in $L^2_{loc}(\R^{N})$.
Then we let $\vp\in C^\infty_c(B_1)$ such that $\vp\equiv 1$ on $B_{\frac{1}{2}}$. By \eqref{eq:Levi}, multiplying  \eqref{eq:weps-stf} by $\vp^2 v_n$  and integrating by parts, we obtain
\begin{eqnarray*}
\displaystyle\int_{ \R^{N+1}_+}|\n (\vp v_n)|^2dx
 &=&
\displaystyle S_{N, s}\int_{ \R^N}|y|^{-s} |v_n|^{q_s-2}|\vp v_n|^2d ydt+o_n(1)\\
&\leq&\displaystyle\frac{S_{N, s}}{2^{\frac{q_s-2}{q_s}}}
\displaystyle\left(   \int_{ \R^N}|y|^{-s} |\vp v_n|^{q_s}dy dt  \right)^{ \frac{2}{q_s}}+o_n(1).
\end{eqnarray*}
Therefore
$$
S_{N, s} \left(   \int_{ \R^N}|y|^{-s} |\vp v_n|^{q_s}dy dt  \right)^{ \frac{2}{q_s}} \leq
\frac{S_{N, s}}{2^{\frac{q_s-2}{q_s}}}
\left(   \int_{ \R^N}|y|^{-s} |\vp v_n|^{q_s}dy dt  \right)^{ \frac{2}{q_s}}+o_n(1).
$$
Making use of the fact that $s\in(0,1)$, we can easily prove that $ S_{N, s}> \frac{S_{N, s}}{2^{\frac{q_s-2}{q_s}}}$. Consequently so that
$$
o_n(1)=\int_{ \R^N}|y|^{-s} |\vp v_n|^{q_s}dy dt=\int_{B_1}|y|^{-s} |v_n|^{q_s}dy dt +o_n(1).
$$
This then contradicts \eqref{eq:Levi}. Therefore $v\neq0$ is a minimizer. By standard arguments, we can show that $v^+=\max(v,0)$ is also a minimizer for $S_{N, s}$. Hence we get the desired result by the maximum principle.
\end{proof}
 %
\subsection{Symmetry and decay estimates of solution of the limiting problem}
\begin{theorem}\label{th:sym-dec}
For $N\geq2$, we let $w\in \calD$ positive, solution of
\begin{align}\label{eq:wgrst}
\begin{cases}
\D w=0& \quad\textrm{ in }  \R^{N+1}_+\vspace{3mm}\\
\displaystyle-\frac{\de w}{\de z}=S_{N, s} |y|^{-s}w^{q_s-1}& \quad\textrm{ on } \R^{N}.
  \end{cases}
\end{align}
where $\R^{N+1}_+=\lbrace
x=(t, y, z) \in \R\times \R^{N-1} \times \R_+
\rbrace$
with boundary 
$$
\de \R^{N+1}_+= \R^N:=\lbrace (t, y )\in \R \times \R^{N-1} \rbrace.
$$
Then
\begin{itemize}
\item[1)] $w(x)=w(t,y, z)$ depends only on $|t|$, $|y|$ and $z$.\\
\noindent
\item[2).] There exist $c_1$ and $c_2$ positive such that 
\begin{equation}\label{Estim1}
\frac{c_1}{1+|x|^{N-1}} \leq w(x)\leq \frac{c_2}{1+|x|^{N-1}} \qquad \textrm{ in $\R^{N+1}_+$}.
\end{equation}
\end{itemize}
\end{theorem}
%
\begin{proof}
\begin{itemize}
\item[1)]
For $\lambda >0$, we define
$$
H_\lambda:= \{ x=(t, y, z )\in \R^{N+1}_+ \::\: t> \lambda\}. 
$$
For $x \in H_\lambda$, we define by
$$
x_\lambda=(2\lambda-t, y,z)
$$
its reflection at the hyperplane $\de H_\lambda$.
Set
$$
u_\lambda: \overline {\R^{N+1}_+ \cap H_{\lambda}} \to \R, \qquad
u_\lambda(x)=w(x_\lambda)-w(z).
$$
Then $u_\lambda$ satisfies
\begin{equation}\label{Ret}
\begin{cases}
-\D u_\lambda=0 \qquad &\textrm{ in $H_\lambda$}\\\\
\displaystyle-\frac{\partial  u_\lambda}{\partial z} = S_{N, s} A(x) \frac{u_\lambda}{|y|^s} \quad & \text{on $\R^N \cap H_\lambda$}\\\\
u_\lambda=0 \quad &\text{on $\R^{N+1}_+ \cap \partial H_\lambda$}.
\end{cases}
\end{equation}
Since the function $t \longmapsto t^{q_s-1}$ is convex on $(0, \infty)$, then we have
$$
0 \leq A:=\frac{w_\lambda^{q-1}-w^{q_s-1}}{w_\lambda-w} \leq (q_s-1) u_\lambda^{q_s-2},
$$
where we have set $v_\lambda=\max(w_\lambda; w)$. Next, we multiply the Euler-Lagrange equation \eqref{Ret} by $u_\lambda^-= \min \{u_\lambda,0\}$ and apply the integration by parts formula to get
\begin{align*}
\int_{H_\lambda} &|\nabla u_\lambda^-|^2\,dz = \int_{H_\lambda} \nabla
u_\lambda \nabla u_\lambda^-\,dz= - \int_{\R^N \cap
  H_\lambda} \frac{\partial u_\lambda}{\partial z} u_\lambda^-\,d\sigma(x)\\
&= S_{N, s}\int_{\R^N  \cap
  H_\lambda} A(x) u_\lambda^- \frac{u_\lambda}{|y|^s} d \sigma(x) \nonumber \\
&\leq  (q_s-1)S_{N, s} \int_{\R^N  \cap  H_\lambda \cap \lbrace u_\lambda \leq 0\rbrace} |u_\lambda^-(z)|^2  
|y|^{-s} w^{q-2}(x) \,d\sigma(x).
\end{align*}
By Hölder's inequality, we obtain
\begin{equation}
  \label{eq:1}
\int_{H_\lambda}  |\nabla u_\lambda^-|^2  
 \le  C(\l)
\Bigl(\int_{\R^N  \cap H_\lambda}|y|^{-s}|u_\lambda^-|^{q_s}\,d \sigma(x) \Bigr)^{\frac{2}{q_s}},
\end{equation}
with
$$
C(\l):=(q-1)S_{N, s}\left(
  \int_{\R^N  \cap H_\lambda\cap \lbrace u_\lambda \leq 0\rbrace}|y|^{-s}w^{q}(z)\,d\sigma(x) \right)^{\frac{q_s-2}{q}}.
$$
Since $c(\lambda) \to 0$ as $\lambda \to \infty$, we have $c(\lambda) <
S$. Therefore, for $\lambda$ sufficiently large, we get
$$
u^\lambda_- \equiv 0 \qquad \textrm{ in} \qquad H_\lambda \cap \R^{N+1}_+.
$$
Set
$$
\lambda^*:= \inf \{\lambda >0 \::\: \text{$w(z) \le
 w(z_{\lambda'})$ for all $z \in H_{\lambda'} \cap \R^{N+1}_+$ and
all $\lambda' \ge \lambda$}\}.
$$
Then $\lambda^*=0$. If $\lambda^*>0$, then $u_{\lambda^*}$ satisfies
$$
\begin{cases}
-\D u_{\lambda^*}=0 \qquad &\textrm{in $\R^{N+1}_+ \cap H_{\lambda^*}$}\\\\
-\displaystyle \frac{\de u_{\lambda^*}}{\de z}=   \frac{w^{q_s-1}(z_{\lambda^*})-w^{q-1}(z)}{|y|^s}
  &\text{on $\R^N \cap H_{\lambda^*}$,}\\\\
u^{\lambda^*} =0   &\textrm{ on 
$\R^{N+1}_+ \cap \partial H_\lambda$}.
\end{cases}
$$
where $\displaystyle \frac{\de u_{\l^*}}{\de z}$ is negative whenever
$w(z_{\lambda*})>0$. 
It follows that, unless $w \equiv 0$, $u^{\lambda^*}$ is strictly positive in $\mathbb{R}^{N+1}_+ \cap H_{\lambda^*}$ by the strong maximum principle. We then pick $D$ sufficiently large so that $D \Subset \mathbb{R}^N \cap H_{\lambda^*}$ and
$$
(q_s-1)S_{N, s}\left(
  \int_{\R^{N} \cap H^{\lambda^*} \setminus D}|z|^{-s}w^{q}(z)\,d\sigma(z) \right)^{\frac{q_s-2}{q_s}} < S_{N, s}.
$$
Then, for $\lambda<\lambda^*$ close to $\lambda^*$,  we have
$D  \subset \R^{N} \cap H_{\lambda}$,
$$
(q_s-1)S_{N, s}\Bigl(
  \int_{\R^{N} \cap H_{\lambda} \setminus
      D}|z|^{-s}w^{q}(z)\,d\sigma(z) \Bigr)^{\frac{q_s-2}{q_s}} < S_{N, s}.
$$
and $u^\lambda >0$ in $D$. As a consequence, $c(\lambda)< S$
for $\lambda<\lambda^*$ close to $\lambda^*$ because $M_\lambda\subset\R^{N} \cap H_{\lambda} \setminus
      D $. By (\ref{eq:1}) we have 
$u^\lambda \ge 0$ in $H_{\lambda} \cap \R^{N+1}_+$ for
$\lambda<\lambda^*$ close to $\lambda^*$,
 contrary to the definition of $\lambda^*$. Therefore
 $$
 w(t, y, -z) \geq w(t, y, z) \qquad \textrm{ in $\R \times \R^{N-1} \times \R_+$}.
 $$
Applying the same argument as before to the function $(t, y, z) \mapsto w(-t, y, z)$, we deduce that
\[
w(-t, y, z) \geq w(t, y, z) \qquad \text{in } \mathbb{R}^{N+1}_+.
\]
Hence 
\begin{equation}\label{Bit1}
w(t, y, z)=w(-t, y, z) \qquad \forall (t, y, z) \in \R^{N+1}_+.
\end{equation}
Repeating the same argument for the functions $y_2 \mapsto w(t, y_2, \cdots, y_N, z)$, with some minor changes, we can easily prove that
$$
w(t, y_2, \cdots, y_{N-1}, z)=w(t, -y_2, \cdots, y_{N-1}, z) \quad \forall (t, y, z) \in \R^{N+1}_+.
$$
and again on the functions $y \longmapsto w(t, By, z)$, where
$B \in O(N)$ is an $(N-1)$-dimensional rotation, we conclude that $w$ only depends on $|t|$, $|y|$ and $z$.
%
\item[2)] For the second point, we write
\begin{align*}
\begin{cases}
\D w =0& \quad \textrm{ in } \R^{N+1}_+\\
- \displaystyle\frac{\de u}{\de z}= a(x) w & \quad \textrm{ on } \R^N,
\end{cases}
\end{align*}
with $a=S_{N, s} |y|^{-s} w^{q-2}\in L^p_{loc}(\R^N)$. Note that, there exists $p>{N}$ such that $a\in L^p_{loc}(\R^N)$. Therefore $w \in L^\infty_{loc} (\R^{N+1}_+)$, see for instance \cite{JLX}. Moreover since the Euler-Lagrange equation \eqref{eq:wgrst} is invariant under Kelvin transform, we immediately  get the desired result. This then ends the proof.
\end{itemize}
%
\end{proof}
\section{Existence of minimizer for $S_{\O, \Sigma, s}$}\label{Section4}
\begin{lemma}\label{lem5}
Let $N \geq 2$. Then for $r>0$ small, there exists $C_r>0$ such that 
\begin{equation}\label{ExtrA}
S_{N, s} \left(\int_{\partial \Omega} \rho_\Sigma^{-s} |u|^{q_s} d\sigma (x)\right)^{2/q_s} \leq (1+r) \int_{\Omega} |\nabla u|^2 dx+ C_r \left[
\int_\Omega u^2 dx+ \left( \int_{\partial \Omega} |u|^{q_s} d \sigma(x)\right)^{2/q_s}
\right].
\end{equation}
\end{lemma}
\begin{proof}
Let $\eta \in \mathcal{C}^1_c(F(Q_{2r}))$ such that 
$$
0 \leq \eta \leq 1 \quad \textrm{ and } \quad \eta \equiv 1 \textrm{ in $F(Q_r)$}.
$$
Since $q_s>2$, we have the existence of some positive constant $C_r>0$ depending on $r$ such that
{\small{\small{\small
$$
\int_{\partial \Omega} \rho_{\Sigma}^{-s} |u|^{q_s} d\sigma(x) \leq (1+r) \left(\int_{\partial \Omega \cap F(Q_r)} \rho_{\Sigma}^{-s} |\eta u|^{q_s} d\sigma(x)\right)^{2/q_s}+C_r \left(\int_{\partial \Omega \cap (F(Q_{2r}) \setminus F(Q_r))} \rho_{\Sigma}^{-s} |(1-\eta)u|^{q_s} d\sigma(x)\right)^{2/q_s}.
$$
}}}
By change of variable formula, Corollary \ref{Det-Metrci-Cor2} and \eqref{Heat1}, we have
\begin{align*}
\int_{\de \O \cap F(Q_r)} \rho_{\Sigma}^{-s} |\eta u|^{q_s} d\sigma(x)&= \int_{\R^{N}} |y|^{-s} |\tilde{\eta} \tilde{u}|^{q_s} \sqrt{|g|}(t, y, 0) dy dt\\\
&\leq (1+cr)\int_{\R^{N}} |y|^{-s} |\tilde{\eta} \tilde{u}|^{q_s} dy dt,
\end{align*}
where we have set
$$
\tilde{\eta}= \eta(F^{-1}(x))\qquad \textrm{and} \qquad \tilde{u}=u(F^{-1}(x)).
$$
Next, thanks to \eqref{HS-Trace}, we obtain
$$
S_{N, s}\left(\int_{\R^{N}} |y|^{-s} |\tilde{\eta} \tilde{u}|^{q_s} dy dt\right)^{2/q_s} \leq \int_{\R^{N+1}_+} |\n (\tilde{\eta} \tilde{u})|^2 dx.
$$ 
Using again a change of variable formula, we can easily see that
$$
\int_{\R^{N+1}_+} |\n (\tilde{\eta} \tilde{u})|^2 dx \leq (1+r) \int_{\Omega} |\n (\eta u)|^2 dx +Cr \int_{\Omega} u^2 dx.
$$
Moreover, we observe that in $\de \O \cap (F(Q_r))$, the function $1-\eta$ vanishes and that $\rho_{\Sigma}$ is bounded in $\de \O \cap (F(Q_{2r}) \setminus F(Q_r))$ so that
$$
\int_{\partial \Omega \cap (F(Q_{2r}) \setminus F(Q_r))} \rho_{\Sigma}^{-s} |(1-\eta)u|^{q_s} d\sigma(x) \leq C_r \int_{\partial \Omega} |u|^{q_s} d\sigma(x),
$$
for some positive constant $C_r$. This then ends the proof.
\end{proof}
\begin{remark}
In general, it's not possible to take the limit as $r \to 0$ in inequality \eqref{ExtrA} except for the Sobolev case $s=0$. We refer to Li and Zhu in  \cite{YanYanLi}: there  exists a positive constant $C=C(\O)$ such that 
\begin{equation}\label{eq:almstST}
\displaystyle  S_{N, 0} \left(\int_{\partial\Omega}| u|^{2^{\sharp}}  d\sigma\right)^{2/2^{\sharp} }
\leq  \int_{\Omega} |\nabla u|^2dx+ C \int_{\partial \Omega}| u|^2d\sigma, \quad \forall u\in H^1(\O).
\end{equation}
\end{remark}
Now, we are in position the prove the following existence result.
\begin{proposition}\label{Existence-Of Solution}
Consider a Lipschitz domain $\Omega \subset \mathbb{R}^{N+1}$ with $N \geq 2$, and let $\Sigma \subset \partial \Omega$ be a closed curve passing through $0$. If $s \in [0,1]$ and $S_{\Omega, \Sigma, s} < S_{N,s}$, then a minimizer for $S_{\Omega, \Sigma, s}$ exists.
\end{proposition}
\begin{proof}
Let $\lbrace u_n\rbrace_{n\geq0} $ be a minimizing sequence for $S_{\Omega, \Sigma, s}$ normalized so that
$$
\int_{\partial \O} \rho_\Sigma^{-s} u_n^{q_s} d\sigma(x)=1 \qquad \textrm{ and } \qquad
 S_{\Omega, \Sigma, s}=\int_\Omega |\nabla u_n|_g^2 dx+\int_\Omega u_n^2 dx+o_n(1).
$$
Then $ \lbrace u_n \rbrace_{n \geq 0}$ is bounded in $ H^1(\O)$ and we assume, up to a subsequence, that
\begin{equation}\label{Convergence}
u_n \rightharpoonup u \quad in\quad H^1(\O); \quad u_n \longrightarrow u \quad in \quad L^{p}(\de \O)  \quad\textrm{ and }\quad u_n \longrightarrow u \quad in \quad L^{p}(\O),
\end{equation}
for  $1 <p < \frac{2N}{N-1}$.
Therefore
\begin{equation}\label{Estimation 1}
S_{\Omega, \Sigma, s} +o_n(1)=\int_\Omega |\nabla u_n| dx+\int_\Omega u^2 dx=\int_\Omega |\nabla u|^2 dx+ \int_\Omega |\nabla (u_n-u)|^2 dx +\int_\Omega u^2 dx +o_n(1).
\end{equation}
By the Brezis-Lieb Lemma \cite{BL}, we have
$$
1=\int_{\partial \Omega} \rho_{\Sigma}^{-s} |u_n|^{q_s} d\sigma(x)= \int_{\partial \Omega} \rho_\Sigma^{-s} |u|^{q_s} d\sigma(x)+\int_{\partial \Omega} \rho_\Sigma^{-s} |u_n-u|^{q_s}  d\sigma(x)+o_n(1).
$$
By Lemma \ref{lem5} and \eqref{Convergence}, we obtain
\begin{equation}\label{Estimation 3}
S_{N,s}\left(\int_{\partial \Omega} \rho_{\Sigma}^{-s} |u_n-u|^{q_s} d\sigma(x)\right)^{2/q_s)} \leq (1+r) \int_\Omega |\nabla (u_n-u)|^2+o_n(1).
\end{equation}
Therefore
$$
S_{N,s} \left(1- \int_{\partial \Omega} \rho_\Sigma^{-s} |u|^{q_s} d\sigma(x)\right)^{2/q_s)}\leq (1+r) \int_\Omega |\nabla (u_n-u)|^2+o_n(1).
$$ 
Using \eqref{Estimation 1} and \eqref{Estimation 3}, we obtain
\begin{equation}
\int_\Omega |\nabla u|^2 dx+ \frac{S_{N,s}}{1+r}\left(1-\int_{\de\Omega} \rho_\Sigma^{-s} |u|^{q_s} d\sigma(x)\right)^{2/q_s}+\int_\Omega u^2 dx \leq S_{\Omega, \Sigma, s}.
\end{equation}
Since
$$
S_{\Omega, \Sigma, s} \biggl( \int_{\partial \Omega}\rho_\Sigma^{-s} |u|^{q_s} d\sigma(x)\biggl)^{2/q_s} \leq \int_\Omega |\nabla u|^2 dx+ \int_\Omega u^2 dx
$$
we get
$$
\frac{S_{N,s}}{1+r}\biggl(1-\int_{\partial \Omega} \rho_\Sigma^{-s} |u|^{q_s} dv_g\biggl)^{2/q_s} \leq S_{\Omega, \Sigma, s} \biggl(1-\biggl( \int_{\de \Omega}\rho_\Sigma^{-s} |u|^{q_s} d\sigma(x)\biggl)^{2/q_s)}\biggl).
$$
Moreover 
$$
\displaystyle 1-\biggl( \int_{\de \Omega}\rho_\Sigma^{-s} |u|^{q_s} dv_g\biggl)^{2/q_s} \leq \left(1-\int_{\de \O} \rho_\Sigma^{-s} |u|^{q_s} d\s(x)\right)^{2/q_s}.
$$
Taking the limit as $r \longrightarrow 0$ we obtain
$$
\bigl( S_{N,\s}- S_{\Omega, \Sigma, s}\bigl) \biggl(1-\biggl( \int_{\de \O}\rho^{-s} |u|^{q_s} dv_g\biggl)^{2/q_s}\biggl) \leq 0.
$$
Since 
$$
S_{N,s}< S_{\Omega, \Sigma, s}\qquad \textrm{ and } \qquad \displaystyle \int_{\de \O}\rho_\Sigma^{-s} |u|^{q_s} d\s(x) \leq 1,
$$
it follows that
$$
\displaystyle \int_{\de \O} \rho_\Sigma^{-s} |u|^{q_s}= 1.
$$
Therefore $u_n \longrightarrow u$ in $H^1(M)$. In particular $u$ is a minimizer for $S_{\Omega, \Sigma, s}$. This then ends the proof.
\end{proof}
\section{Construction of Test Function and Comparing Hardy-Sobolev Best Constants}\label{Section5}
\subsection{Proof of Theorem \ref{Theorem1}}
Let $w \in \calD^{1,2}(\R^{N+1}_+)$ be a positive solution for
\begin{equation}\label{Euler-Lagrange-In-RN}
\begin{cases}
\Delta w=0 \qquad &\textrm{in $\R^{N+1}_+$}\\\
\displaystyle -\frac{\de w}{\de z}= S_{N, s} |y|^{-s} w^{q_s-1} &\textrm{on $\R^{N}$}.
\end{cases}
\end{equation}
For $\e>0$, we define
$$
v_\e(F(x)):=\e^{\frac{1-N}{2}} w(x/\e).
$$
Next, we let $\eta \in \calC^\infty_c(F(\mathcal{Q}_{2r}))$ such that
\begin{equation}\label{Cutt-Off-Function}
\eta \equiv 1 \qquad \textrm{ in $F(\mathcal{Q}_{2r})$ and } 0\leq \eta \leq 1.
\end{equation}
Then we define the test function as
\begin{equation}\label{Test-Function}
u_\e(F(x))=\eta(F(x)) v_\e(F(x)).
\end{equation}
Therefore $u_\e \in H^1(\O)$. In the following, we will expand the functional $J: H^1(\Omega) \to \R$ defined by
$$
J(u)= \frac{\displaystyle \int_{\Omega} |\n u|^2 dx+ \int_{\Omega} u^2 dx}{\left(\displaystyle\int_{\de \O} \rho_{\Sigma}^{-s} |u|^{q_s} d\s(x)\right)^{2/q_s}}.
$$
\begin{lemma}\label{Gradient}
We have
\begin{align}\label{El7}
\int_{\O} |\n u_\e|^2 dx+\int_{\O} u_\e^2 dx=\int_{\R^N} |\n w|^2 dx+\e H_{\de \Omega}(0) &\int_{\mathcal{Q}_{r/\e}} z |\n w|^2 dx-2\e \mathcal{H}_{1} \int_{\mathcal{Q}_{r/\e}} z |\frac{\de w}{\de t}|^2 dx\nonumber\\\
&-2\e\mathcal{H}_{\de \O}(0) \int_{\mathcal{Q}_{r/\e}} z |\n_y w|^2 dx+o(\e),
\end{align}
where the geometric quantities $H_{\de \O}$, $\mathcal{H}_{\de \O}$ and $\mathcal{H}_1$ are defined in \eqref{Mean-Curvature} and \eqref{Rel-Mean-Curvature}. 
\end{lemma}
\begin{proof}
We have
\begin{equation}\label{El1}
\int_\O |\n u_\e|^2 dx= \int_{\O\cap F(\mathcal{Q}_r)} |\n u_\e|^2 dx+\int_{\O \cap (F(\mathcal{Q}_{2r} \setminus F(\mathcal{Q}_r)} |\n u_\e|^2 dx.
\end{equation}
By the change of variable formula $\tilde{x}=\frac{F^{-1}(x)}{\e}$ and \eqref{Cutt-Off-Function}, we have
\begin{align}\label{El2}
\int_{\O\cap F(\mathcal{Q}_r)} & |\n u_\e|^2 dx=\sum_{ij=1}^{N+1} \int_{\mathcal{Q}_{r/\e}} g^{ij}(\e x) \frac{\de w}{\de x_i} \frac{\de w}{\de x_j} \sqrt{|g|}(\e x) dx\nonumber\\\
&=\int_{\mathcal{Q}_{r/\e}} |\n w|^2 \sqrt{|g|}(\e x) dx+\sum_{ij=1}^{N+1} \int_{\mathcal{Q}_{r/\e}} \left(g^{ij}(\e x)-\delta_{ij} \right) \frac{\de w}{\de x_i} \frac{\de w}{\de x_j} \sqrt{|g|}(\e x) dx.
\end{align} 
By Corollary \ref{Det-Metrci-Cor2} and using the fact that $w$ depends on $|y|$, we get
\begin{align}\label{El3}
&\int_{\mathcal{Q}_{r/\e}} |\n w|^2 \sqrt{|g|}(\e x) dx=\int_{\mathcal{Q}_{r/\e}} |\n w|^2 dx+\e H_{\de \Omega}(0) \int_{\mathcal{Q}_{r/\e}} z |\n w|^2 dx +O\left(\e^2 \int_{\mathcal{Q}_{r/\e}} |x|^2 |\n w|^2 dx\right)\nonumber\\\
&=\int_{\R^N} |\n w|^2 dx+\e H_{\de \Omega}(0) \int_{\mathcal{Q}_{r/\e}} z |\n w|^2 dx +O\left(\e^2 \int_{\mathcal{Q}_{r/\e}} |x|^2 |\n w|^2 dx+ \int_{\R^N \setminus \mathcal{Q}_{r/\e}} |\n w|^2 dx\right).
\end{align}
We have
\begin{align*}
&\sum_{ij=1}^{N+1} \int_{\mathcal{Q}_{r/\e}} \left(g^{ij}(\e x)-\delta_{ij} \right) \frac{\de w}{\de x_i} \frac{\de w}{\de x_j} \sqrt{|g|}(\e x) dx= \int_{\mathcal{Q}_{r/\e}} \left(g^{11}(\e x)-1 \right) |\frac{\de w}{\de t}|^2 \sqrt{|g|}(\e x) dx\\\\
&+ \int_{\mathcal{Q}_{r/\e}} \left(g^{N+1 N+1}(\e x)-1 \right) |\frac{\de w}{\de z}|^2 \sqrt{|g|}(\e x) dx+\sum_{ij=2}^{N}\int_{\mathcal{Q}_{r/\e}} \left(g^{ij}(\e x)-\delta_{ij} \right) |\n_y w|^2 \frac{y_i y_j}{|y|^2} \sqrt{|g|}(\e x) dx \\\\
&\displaystyle+\sum_{i=2}^{N} \int_{\mathcal{Q}_{r/\e}} g^{1i}(\e x)\n_y w\frac{y_i}{|y|} \frac{\de w}{\de t}  \frac{t}{|t|}\sqrt{|g|}(\e x) dx+\sum_{i=2}^{N} \int_{\mathcal{Q}_{r/\e}} g^{N+1 i}(\e x)\n_y w\frac{y_i}{|y|} \frac{\de w}{\de z} \sqrt{|g|}(\e x) dx\\\\
&\qquad \qquad \qquad \qquad \qquad \qquad \qquad \qquad \qquad \qquad+ 
 \int_{\mathcal{Q}_{r/\e}} g^{N+1 1}(\e x)\frac{\de w}{\de z} \frac{\de w}{\de t}  \frac{t}{|t|}\sqrt{|g|}(\e x) dx.
\end{align*}
Therefore by Corollary \ref{Inverse-Metric-Cor1},  Corollary \ref{Det-Metrci-Cor2} and using the symmetry properties of $w$ given by Theorem \ref{th:sym-dec}, we obtain
\begin{align}\label{El4}
\sum_{ij=1}^{N+1} \int_{\mathcal{Q}_{r/\e}} \left(g^{ij}(\e x)-\d_{ij} \right) \frac{\de w}{\de x_i} \frac{\de w}{\de x_j}& \sqrt{|g|}(\e x) dx=-2\e H_{11} \int_{\mathcal{Q}_{r/\e}} z |\frac{\de w}{\de t}|^2 dx\nonumber\\\
&-\frac{2}{N-1}\e \sum_{i=2}^N H_{ii}(0) \int_{\mathcal{Q}_{r/\e}} z |\n_y w|^2 dx+O\left(\e^2 \int_{\mathcal{Q}_{r/\e}} |x|^2 |\n w|^2 dx \right).
\end{align}
Combining \eqref{El2}, \eqref{El3} and \eqref{El4} we obtain
\begin{align}\label{El5}
\int_{\O\cap F(\mathcal{Q}_r)} |\n u_\e|^2 dx=\int_{\R_+^{N+1}} |\n w|^2 dx+\e H_{\de \Omega}(0) &\int_{\mathcal{Q}_{r/\e}} z |\n w|^2 dx-\e H_{11} \int_{\mathcal{Q}_{r/\e}} z |\frac{\de w}{\de t}|^2 dx\nonumber\\\
&-\frac{2}{N-1}\e \sum_{i=2}^N H_{ii}(0) \int_{\mathcal{Q}_{r/\e}} z |\n_y w|^2 dx+O\left(\rho(\e)\right),
\end{align}
where
\begin{equation}\label{Max0}
\rho(\e)=\e^2 \int_{\mathcal{Q}_{r/\e}} |x|^2 |\n w|^2 dx+ \int_{\R^{N+1}_+ \setminus \mathcal{Q}_{r/\e}} |\n w|^2 dx.
\end{equation}
By the change of variable formula $\tilde{x}=\frac{F^{-1}(x)}{\e}$, we have
\begin{align}\label{El6}
\int_{\O \cap (F(\mathcal{Q}_{2r} \setminus F(\mathcal{Q}_r)} |\n u_\e|^2 dx= \int_{\mathcal{Q}_{2r/\e} \setminus \mathcal{Q}_{r/\e}} |\n w|^2 dx.
\end{align}
Combining \eqref{El1}, \eqref{El5} and \eqref{El6}, we obtain
\begin{align}\label{El7}
\int_{\O\cap F(\mathcal{Q}_r)} |\n u_\e|^2 dx=\int_{\R_+^{N+1}} |\n w|^2 dx+\e H_{\de \Omega}(0) &\int_{\mathcal{Q}_{r/\e}} z |\n w|^2 dx-2\e H_{11} \int_{\mathcal{Q}_{r/\e}} z |\frac{\de w}{\de t}|^2 dx\nonumber\\\
&-\frac{2}{N-1}\e \sum_{i=2}^N H_{ii}(0) \int_{\mathcal{Q}_{r/\e}} z |\n_y w|^2 dx+O\left(\rho(\e)\right).
\end{align}
In the sequel, we will estimate the error term $\rho(\e)$. For that we let $\eta \in \calC^\infty_c(Q_{2r})$ radial such that $\eta\equiv1$ in $Q_r$ and $0\leq \eta \leq 1$. We define
$$
\eta_\e(x)= \eta(\e x).
$$
Then we multiply \eqref{Euler-Lagrange-In-RN} by $|x|^2 \eta_\e w$ and apply the integration by parts formula to get
$$
0=\int_{Q_{2r/\e}} \eta_\e |x|^2 w \Delta w dx=-\int_{Q_{2r/\e}} \n w \cdot \n (\eta_\e |x|^2 w) dx +\int_{\R^N \cap Q_{2r/\e}} \frac{\de w}{\de \nu} \left( \eta_\e |(t, y)|^2 w\right) dy dt.
$$
Therefore
\begin{align*}
\int_{Q_{r/\e}}|x|^2 |\n w|^2 dx&=  O\left(\int_{Q_{2r/\e}} \n w^2 \cdot \n (\eta_\e |x|^2) dx+\int_{\R^N \cap Q_{2r/\e}} \eta_\e |y|^{2-s} w^{q_s} dy dt\right)\\\\
&=O\left(\int_{Q_{r/\e}} w^2dx+\int_{Q_{2r/\e} \setminus Q_{r/\e}} w^2 \D(\eta_\e |x|^2)dx+\int_{\R^N \cap Q_{2r/\e}} |y|^{2-s} w^{q_s} dy dt\right)\\\\
&=O\left(\int_{Q_{r/\e}} w^2dx+\int_{Q_{2r/\e} \setminus Q_{r/\e}} w^2 (1+\e^2 |x|^2+\e)dx+\int_{\R^N \cap Q_{2r/\e}} |y|^{2-s} w^{q_s} dy dt\right).
\end{align*}
Hence thanks to the estimate \eqref{Estim1} and using polar coordinates, it easy follows that
\begin{equation}\label{Max1}
\int_{Q_{r/\e}}|x|^2 |\n w|^2 dx \sim C + 
\begin{cases}
O(\e^{N-3}) \qquad &\textrm{ if $N \geq 4$}\\\
O(ln(\e)) &\textrm{ if $N=3$}.
\end{cases}
\end{equation}
Next, we multiply \eqref{Euler-Lagrange-In-RN} by $(1-\eta_\e) w$ and apply the integration by parts formula to get
\begin{align*}
0=\int_{\R^{N+1}_+} (1-\eta_\e) w \D w dx= \int_{\R^{N+1}_+} (1-\eta_\e) |\n w|^2 dx &+\frac{1}{2} \int_{\R^{N+1}_+} \n
 w^2 \cdot \n(1-\eta_\e) dx\\\
 & = \int_{\R^{N}} |y|^{-s} w^{q_s} (1-\eta_\e) dy dt.
\end{align*}
 Therefore
\begin{align*}
\int_{\R^{N+1}_+\setminus Q_{r/\e}} |\n w|^2 dx&=O\left(\int_{\setminus Q_{2r/\e} \setminus Q_{r/\e}}\D \eta_\e w^2 dx+ \int_{\R^{N} \setminus \left((-\frac{r}{\e}, \frac{r}{\e}) \times B^{N-1}_{r/\e}\right)} |y|^{-s} w^{q_s} dydt\right)\\\
&=O\left(\e^2 \int_{Q_{2r/\e} \setminus Q_{r/\e}} w^2 dx+ \int_{\R^{N}\setminus B^{N}_{r/\e}} |y|^{-s} w^{q_s} dydt\right).
\end{align*}
Hence thanks to the estimate  \eqref{Estim1} and using polar coordinates, it easy follows that
$$
\int_{Q_{2r/\e} \setminus Q_{r/\e}} w^2 dx= O(1) \qquad  \textrm{ and } \qquad \int_{\R^{N}\setminus B^{N}_{r/\e}} |y|^{-s} w^{q_s} dydt=O(\e^{N-s}) \qquad  \forall N \geq 3.
$$
Therefore
\begin{equation}\label{Max2}
\int_{\R^{N+1}_+\setminus Q_{r/\e}} |\n w|^2 dx=O(\e^2) \qquad \textrm{ for all $N \geq 3$}.
\end{equation}
By \eqref{Max0}, \eqref{Max1} and \eqref{Max2}, we obtain, for all $N \geq 3$, that
\begin{equation}\label{Ndongo1}
\rho(\e)= o(\e), \qquad \textrm{ as $\e \to 0$}.
\end{equation}
We finish the proof by estimating $\int_{\O} u_\e^2 dx$. By change of variable formula, we have
$$
\int_{\O} u_\e^2 dx= O\left( \e^2 \int_{\mathcal{Q}_{r/\e}} w^2 dx\right)
$$
and immediately from \eqref{Estim1}, we obtain
$$
\int_{\mathcal{Q}_{r/\e}} w^2 dx\sim C+
\begin{cases}
O(\ln(\e)) \qquad &\textrm{ if $N=3$}\\\
O(\e^{N-3}) & \textrm{ if $N \geq 4$}.
\end{cases}
$$
Therefore
\begin{equation}\label{Ndongo2}
\int_{\O} u_\e^2 dx=o(\e) \qquad \textrm{ for all $N \geq 3$}.
\end{equation}
Hence by \eqref{El1}, \eqref{El5}, \eqref{El6}, \eqref{El7}, \eqref{Ndongo1} and \eqref{Ndongo2}, we obtain the desired result. This then ends the proof.
\end{proof}

\begin{lemma}\label{Denominator}
We have
$$
\int_{\de \O} \rho_{\Sigma}^{-s} |u_\e|^{q_s} d\sigma(x)= \int_{\R \times \R^{N-1}} |y|^{-s} |w|^{q_s} dt dy+o(\e).
$$
\end{lemma}
\begin{proof}
We have
$$
\int_{\de \Omega} \rho^{-s}_{\Sigma} |u_\e|^{q_s} d\s (x)= \int_{\de \Omega \cap F(\mathcal{Q}_r)} \rho^{-s}_{\Sigma} |u_\e|^{q_s} d\s (x)+\int_{\de \Omega \cap F(\mathcal{Q}_{2r}) \setminus F(\mathcal{Q}_{r}) )} \rho^{-s}_{\Sigma} |u_\e|^{q_s} d\s (x).
$$
By the change of variable formula $\tilde{x}=\frac{F^{-1}(x)}{\e}$, \eqref{Heat1} and Corollary \ref{Det-Metrci-Cor2}, we have
\begin{align}
\int_{\de \Omega \cap F(\mathcal{Q}_r)} \rho^{-s}_{\Sigma} &|u_\e|^{q_s} d\s (x)= \int_{-r/\e}^{-r/\e}\int_{B(0, r/\e)} |y|^{-s} |w(t, y, 0)|^{q_s} \sqrt{|g|}(\e t, \e y, 0) dt dy\nonumber\\\
&=\int_{-r/\e}^{-r/\e}\int_{B(0, r/\e)} |y|^{-s} |w(t, y, 0)|^{q_s} \left(1+\e \sum_{i=2}^N \kappa_i y_i +O(\e^2 |(t, y)|^2)\right) dt dy\nonumber\\\
&=\int_{\R^N} |y|^{-s} |w(t, y, 0)|^{q_s} dt dy+O\left(\rho_2(\e)\right),
\end{align}
where
$$
\rho_2(\e):=\e^2 \int_{\mathcal{B}_{r/\e}} |(t,y)|^{2} |y|^{-s } |w(t, y, 0)|^{q_s} dt dy+ \int_{\R^N \setminus \mathcal{B}_{r/\e}} |y|^{-s} |w(t, y, 0)|^{q_s} dt dy,
$$
with $\mathcal{B}_{R}:=(-R, R)\times B_{\R^{N-1}}(0, R)$. By the estimate in Theorem \ref{th:sym-dec}, we get
$$
\rho_2(\e)=o(\e) \qquad \textrm{ as $\e \to 0$}.
$$
By change of variable formula and  Theorem \ref{th:sym-dec}, we then have
$$
\int_{\partial \Omega \cap F(\mathcal{Q}_{2r}) \setminus F(\mathcal{Q}_{r}) )} \rho^{-s}_{\Sigma} |u_\epsilon|^{q_s} d\sigma (x)= O \left(\int_{\B_{2r/\epsilon} \setminus \B_{r/\e}} |y|^{-s} w^{q_s} dy dt \right)=o(\e) \qquad \textrm{ as $\e \to 0$}.
$$
This then ends the proof.

\end{proof}
Next, we consider the following Hardy-Sobolev trace constants $S_{N, s}$ and $S_{\Omega, \Sigma, s}$ defined respectively by
$$
S_{N, s}:= \inf \left\lbrace \int_{\R^{N+1}_+} |\n u|^2 dx \qquad \textrm{ such that 
$u \in \calD^{1,2}(\R^{N+1}_+)$ and $\int_{\R^{N}} |y|^{-s} w^{q_s} dy dt=1$}
\right\rbrace.
$$
and
$$
S_{\O, \Sigma,s}=\inf_{u \in H^1(\Omega) \setminus \lbrace 0\rbrace} \frac{\displaystyle\int_\O |\n u|^2 dx+\int_\O u^2 dx}{\left(\displaystyle \int_{\de \O} \rho_{\Sigma}^{-s} |u|^{q_s} d\s(x)\right)^{2/q_s}}.
$$
We conclude this section with the following result, which follows directly from Lemmas~\ref{Gradient} and~\ref{Denominator}.
\begin{proposition}\label{PropreAZER}
Let $N \geq 3.$ Then we have
\begin{align*}
\mathcal{S}_{\O,s}^\Sigma \leq  S_{N, s}&+\e H_{\de \O}(0) \int_{\R^{N+1}_+} z |\n w|^2 dx\\\\
&-2\e\mathcal{H}_{\de \Omega}(0) \int_{\R^{N+1}_+} z |\n_y w|^2 dx-\e \mathcal{H}_1 \int_{\R^{N+1}_+} z |\frac{\de w}{\de t}|^2 dx+o(\e) \quad \textrm{ as $\e \to 0$},
\end{align*}
where $H_{\de \O}(0)$ is the mean curvature of $\de \Omega$ defined in \eqref{Mean-Curvature},
$$
\mathcal{H}_{\de \Omega}:=\displaystyle\frac{\displaystyle\sum_{i=2}^N\langle H(X_i), X_i \rangle(0)}{N-1} \quad \textrm{ and } \quad  \mathcal{H}_1(0):=\langle H(X_1), X_1\rangle(0).
$$
\end{proposition}
\begin{proof}
Let $\eta \in \calC^\infty_c(Q_{2r})$ radial such that $\eta\equiv1$ in $Q_r$ and $0\leq \eta \leq 1$. We define
$$
\eta_\e(x)= \eta(\e x).
$$
Then we multiply \eqref{Euler-Lagrange-In-RN} by $z \eta_\e w$ and apply the integration by parts formula to get
$$
0=\int_{Q_{\frac{2r}{\e}}} \eta_\e z w \D w dx=-\int_{Q_{\frac{2r}{\e}}}  \n w \cdot \n (z \eta_\e w) dx
$$
which implies that
\begin{align*}
\int_{Q_{\frac{2r}{\e}}} \eta_\e z |\n w|^2 dx&= \frac{1}{2} \int_{Q_{\frac{2r}{\e}}} \n w^2 \cdot \n(z \eta_\e) dx\\\\
&= -\int_{Q_{\frac{2r}{\e}}} w^2 \D(z \eta_\e) dx+ \int_{-r/\e}^{r/\e} \int_{B^{N-1}(0, r/\e)} w^2(t, y, 0) \eta_\e(t, y, 0) dy dt.
\end{align*}
Therefore
\begin{align*}
\int_{Q_{\frac{r}{\e}}} z |\n w|^2 dx=O\left(\int_{Q_{\frac{2r}{\e}}} w^2 \left(z \e^2+ \e|\n \eta|\right) dx+ \int_{B^{N}(0, r/\e)} w^2(t, y, 0) dy dt \right).
\end{align*}
Using \eqref{Estim1} and polar coordinates as previously, we obtain
$$
\int_{Q_{\frac{r}{\e}}} z |\n w|^2 dx \sim C+ 
\begin{cases}
O(\e ln(\e)) \qquad &\textrm{ if $N=3$}\\\
O(\e^{N}) & \textrm{ if $N \geq 4$}.
\end{cases}, \qquad \textrm{ as $\e \to 0$}.
$$
Consequently
$$
\int_{\R^{n+1}_+} z |\n w|^2 dx <\infty \qquad \textrm{ for all $N \geq 3$}.
$$
Hence
$$
\mathcal{S}_{\O,s}^\Sigma=S_{N, s}-\e\mathcal{H}_{\de \Omega}(0) \int_{\R^{N+1}_+} z |\n w|^2 dx-\e \mathcal{H}_1 \int_{\R^{N+1}_+} z |\frac{\de w}{\de t}|^2 dx+o(\e).
$$
\end{proof}
\begin{proof}\textbf{of Theorem \ref{Theorem1}}
Let $(u_\e)_\e \subset H^1(\O)$ defined by \eqref{Test-Function}. Then by definition, we have
$$
\mu_{\Omega, \Sigma, s} \leq J(u_\e)= 
\frac{\displaystyle\int_{\O} \left(|\n u_\e|^2+ u_\e^2\right) dx}{\left(\displaystyle\int_{\de \O} \rho^{-s}_{\Sigma} |u_\e|^{q_s} d\s(x)\right)}.
$$
Moreover by Proposition \ref{PropreAZER}, assuming that
$$
 H_{\de \O}(0)-2\mathcal{H}_{\de \Omega}(0) \frac{\displaystyle\int_{\R^{N+1}_+} z |\n w|^2 dx}{\displaystyle\int_{\R^{N+1}_+} z |\n w|^2 dx}-\mathcal{H}_1 \frac{\displaystyle\int_{\R^{N+1}_+} z |\frac{\de w}{\de t}|^2 dx}{\displaystyle\int_{\R^{N+1}_+} z |\n w|^2 dx}<0,
$$
we deduce that
$$
S_{\Omega, \Sigma, s} <S_{N, s}.
$$
Hence by Proposition \ref{Existence-Of Solution}, there exists a positive function $u \in H^1(\O)$ such that
\begin{align*}
	\begin{cases}
	\displaystyle-\D u+u=0&\qquad \textrm{ in } \O \vspace{3mm}\\
	\displaystyle \frac{\de u}{\de \nu}=S_{\Omega, \Sigma, s}\, \rho_\Sigma^{-s}(\s)\,  u^{q_s-1}& \qquad \textrm{ on } \de \O,
	\end{cases}
	\end{align*}
	where $\nu$ is the unit outer normal of $\de\O$. This then ends the proof.
\end{proof}
\subsection{Proof of Theorem \ref{Theorem2}}.

\begin{proof}
The proof is Theorem \ref{Theorem2} is very simple. Indeed, since the domain $\O$ is bounded, the constants functions are in $H^1(\O)$. Recall that 
$$
J(u)=\frac{\displaystyle\int_{\O} \left(|\n u|^2+ u^2\right) dx}{\left(\displaystyle\int_{\de \O} \rho^{-s}_{\Sigma} |u|^{q_s} d\s(x)\right)}.
$$
Then for any $c\in \R^*$, we have
$$
S_{\Omega, \Sigma, s} \leq J(c)=\frac{|\O|}{\left(\displaystyle  \int_{\Omega} \rho^{-s}_{\Sigma} d\s(x)\right)^{2/q_s}}.
$$
Therefore by  Proposition \ref{Existence-Of Solution}, the constant $S_{\O, \Sigma, s}$ is achieved if 
$$
\frac{|\O|}{\left(\displaystyle  \int_{\Omega} \rho^{-s}_{\Sigma} d\s(x)\right)^{2/q_s}}< S_{N, s}.
$$
This then ends the proof.
\end{proof}
%
%
%
%
%
%
%
%
%
%
%
%
%
%
%
%
%
%
%
%
%
%
%
%
%
%
%
%
%
%
%
%
%
%
%
%

%
%

\section*{Acknowledgments}
The authors thank the anonymous referee for their careful and valuable feedback. Special thanks are also due to Professor Abdoul Salam Diallo for taking the time to read, revise, and help improve the manuscript.

\end{document}